\documentclass[a4paper]{amsart}

\usepackage{amscd}
\usepackage{amssymb}
\usepackage{amsthm}

\title[Universal formulas on Hilbert schemes]{Universal formulas for
  characteristic classes on the Hilbert schemes of points on surfaces}
\author{Samuel Boissi\`ere \and Marc A.~Nieper-Wi\ss kirchen}
\date{\today}
\address{Institut f\"ur Mathematik\\Johannes-Gutenberg Universit\"at\\
  55099 Mainz\\Germany}
\email{boissiere@mathematik.uni-mainz.de\\
  nieper@mathematik.uni-mainz.de}

\theoremstyle{plain}
\newtheorem{theorem}{Theorem}
\newtheorem{lemma}{Lemma}
\newtheorem{corollary}{Corollary}
\newtheorem{proposition}{Proposition}
\theoremstyle{definition}
\newtheorem{definition}{Definition}
\theoremstyle{remark}
\newtheorem{remark}{Remark}
\newtheorem{example}{Example}

\newcommand{\abs}[1]{|#1|}
\newcommand{\norm}[1]{\|#1\|}
\newcommand{\colim}{\varinjlim}
\newcommand{\ch}{\mathrm{ch}}
\newcommand{\td}{\mathrm{td}}

\newcommand{\set}[1]{\mathbf #1}
\newcommand{\sheaf}[1]{\mathcal #1}
\newcommand{\Grass}{\set G}
\newcommand{\SG}{\mathfrak S}

\newcommand{\one}{\mathbf{\left|1\right>}}
\newcommand{\vac}{\mathbf{\left|0\right>}}
\newcommand{\Hilb}[2]{#2^{[#1]}}
\newcommand{\Hilbtot}{\operatorname{Hilb}}
\newcommand{\End}{\operatorname{End}}
\newcommand{\supp}{\operatorname{supp}}
\newcommand{\id}{\mathrm{id}}
\newcommand{\Shi}{\operatorname{Shi}}

\begin{document}

\begin{abstract}
  This article can be seen as a sequel to the first author's article
  ``Chern classes of the tangent bundle on the Hilbert scheme of
  points on the affine plane'', where he calculates the total Chern
  class of the Hilbert schemes of points on the
  affine plane by proving a result on the existence of certain
  universal formulas expressing characteristic classes on the Hilbert
  schemes in term of Nakajima's creation operators.
  
  The purpose of this work is (at least) two-fold. First of all, we
  clarify the notion of ``universality'' of certain formulas about the
  cohomology of the Hilbert schemes by defining a universal algebra of
  creation operators. This helps us to reformulate and extend a lot of
  the first author's previous results in a very precise manner.

  Secondly, we are able to extend the previously found results by
  showing how to calculate any characteristic class of the Hilbert
  scheme of points on the affine plane in terms of the creation
  operators. In particular, we have included the calculation of the
  total Segre class and the square root of the Todd class.
  
  Using this methods, we have also found a way to calculate any
  characteristic class of any tautological sheaf on the Hilbert scheme
  of points on the affine plane. This in fact gives another complete
  description of the ring structure of the cohomology spaces of the
  Hilbert schemes of points on the affine plane.
\end{abstract}



\maketitle

\tableofcontents

\section{Introduction}

The Hilbert scheme of $n$ points on a complex surface $X$
(see~\cite{NHS} or~\cite{CM} for the non-algebraic case) parametrises
the zero-dimensional subspaces of length $n$ on $X$. It is denoted in
this text by $X^{[n]}$. By a result of Fogarty (\cite{FOG}), it is
smooth of dimension $2n$. The vector space structure of its rational
cohomology ring has been calculated by G\"ottsche (\cite{GOE}),
i.e.~he calculated the Betti numbers. His formulas can be a given a
particularly nice form if one considers all Hilbert schemes at once,
i.e.~if one studies $\Hilbtot X := \amalg_{n \ge 0} \Hilb n X$.
Nakajima (\cite{NAK}) and Grojnowski (\cite{GRO}) were led by this
result to the construction of certain operators $q_n(\alpha)$ and
$p_n(\alpha)$ on the cohomology $\set H_{\set Q} X$ of $\Hilbtot X$
that are part of a vertex algebra structure (see~\cite{LEH}). Here
$\alpha$ is a cohomology class on the surface $X$ and $n$ is a
positive integer. The upshot of this construction is that every
cohomology class in $\set H_{\set Q} X$ can be given as an application
of a polynomial in the $q_n(\alpha)$ on a certain element $\vac$ in
$\set H_{\set Q} X$, called the vacuum. In the case of the affine
plane, i.e.~$X = \set C^2$, this leads to the isomorphism $\set
H_{\set Q} \set C^2 \simeq \set Q[[p_1, p_2, \dots]]$, the power
series ring in infinitely many variables, as vector spaces (not rings!).

Each sheaf $\sheaf F$ on $X$ gives rise to a sheaf $\Hilb n {\sheaf
  F}$ on each $\Hilb n X$, called a tautological sheaf. One can ask
what the total Chern classes of these sheaves in terms of polynomials
in the $q_n(\alpha)$ applied to the vacuum are. Lehn gives
in~\cite{LEH} a closed formula for the case that $\sheaf F$ is an
invertible sheaf.

Another natural sheaf on the Hilbert schemes is given by the tangent
sheaf. The first author has started to attack the problem of
expressing the total Chern classes of these tangent sheaves in terms
of the $q_n(\alpha)$ in~\cite{BOS2}. It turns out that these
expressions do not really depend on the surface $X$ (or the invertible
sheaf $\sheaf F$) but only on some algebraic invariants of $X$ and $\sheaf F$
(namely more or less their Chern classes). Thus one can say there are
\emph{universal formulas} governing the Chern classes of the
tautological and the tangent sheaves. Unfortunately the coefficients
in these universal formulas are not known in a closed form for the
case of the tangent sheaves of the Hilbert schemes on an arbitrary
surface (or the tautological sheaves of a sheaf that is not of rank
one). However, when one specialises to $X = \set C^2$, the affine
plane, there is a closed formula for the total Chern class of the
tangent sheaves on the Hilbert schemes over $X$ (expressed as elements
in the ring $\set Q[[p_1, p_2, \dots]]$), which is the main
result of~\cite{BOS2}.

This article builts upon these results. We clarify what is meant by
universality of certain formulas by introducing a universal algebra of
operators together with an evaluation map for each surface $X$ that
specialises these operators to the operators $q_n(\alpha)$. Equipped
with this machinery, we easily rediscover the results of~\cite{BOS2}
that state how certain formulas have to look like.

We improve the methods that were used in~\cite{BOS2} to
calculate the total Chern classes of the tangent sheaves
$\Theta_{\Hilb n X}$ of the Hilbert schemes of points on the affine
plane to be able to calculate all multiplicative classes of
$\Theta_{\Hilb n X}$, i.e.~the total Segre class or Todd class. Given
a multiplicative class $\gamma = \prod_{i = 1}^\infty f(\lambda_i)$
where the $\lambda_i$ are the (universal) Chern roots and $f \in 1 + x
R[[x]]$ is a power series with coefficients in a ring $R$, one of our
main results states that
$$\sum_{n \ge 0} \gamma(\Theta_{\Hilb n{(\set C^2)}}) = \exp\left(
  \sum_{k \ge 1} g_{2k + 1}\cdot\frac{p_{2k + 1}}{2k + 1}\right) \in
R[[p_1, p_2, \dots]],$$
where
the power series $g = \sum_{n \ge 1} g_n t^n$ is defined by
$\frac{\partial g}{\partial t}\left(\frac x{f(x) f(-x)}\right) =
f(x)f(-x)$. As every characteristic class appears in a sufficiently general
multiplicative class, one can use this result to calculate effectively
all characteristic classes of the sheaves $\Theta_{\Hilb n {(\set
    C^2)}}$.

As mentioned before, Lehn has calculated the total Chern classes of
the tautological sheaves in the rank one case. We are able to extend
his result to arbitrary multiplicative classes in the case of the
affine plane (where there are only the tautological sheaves of the
structure sheaf). Our second main result states that
$$\sum_{n \ge 0} \gamma(\Hilb n {\sheaf O}) = \exp\left( \sum_{n \ge
    1} g_n\cdot\frac{p_n}{n}\right) \in R[[p_1, p_2, \dots]],$$
where
this the power series $g = \sum_{n \ge 1} g_n t^n$ is defined by
$\frac{\partial g}{\partial t}\left(\frac x{f(-x)}\right) = f(-x)$. As
explained in the last section of this article, one can use this result
to effectively calculate the cup-product in $\set H_R \set C^2 \simeq
R[[p_1, p_2, \dots]]$.

It is remarkable that these formulas for the case of the affine plane
do also give partial results for the characteristic classes of the
tautological sheaves and the tangent sheaves of the Hilbert schemes of
points on any surface by virtue of the universality and uniqueness of
the involved formulas. This is explained in the main text.

Let us outline the structure of this article. The first of the
sections that follow is dedicated to the introduction of some
notations and notions that are used throughout the text.

The third section gives some general results on (multiplicative)
characteristic classes. As we want to discuss characteristic classes
of the tangent sheaf of a complex manifold and characteristic classes
of bundles on it uniformly, we introduce the notion of a marked
complex manifold $(X, \alpha)$, where $X$ is a complex manifold and
$\alpha$ is a cohomology class on $X$. In the applications, $\alpha$
will often by the Chern character of a sheaf on $X$, from which all
other characteristic classes can be built. We introduce the ring of
all (virtual) marked complex manifolds.

In the fourth section, we come to the Hilbert scheme of points on a
surface. We give the definition and the construction of one half of
Nakajima's operators (the half that interests us). It ends with
defining the notion of a Hilbert scheme of a marked complex surface
and explaining that the Hilbert scheme operation is like an
exponential mapping from the abelian group of (virtual) marked complex
surfaces to the ring of all (virtual) marked complex manifolds.

The fifth section is devoted to the universality of certain
formulas, which has already been stated in~\cite{BOS2}. Here, we add
the remark that the universal expressions for multiplicative classes
evaluated on the Hilbert schemes are unique, which makes it possible
to get results for the general case from specific results on the
Hilbert scheme of points on the affine plane.

The last section contains the two main results in that have been
stated above. An algebraic identity that is used in the proof of these
two theorems (and that makes it possible to generalise the result
of~\cite{BOS2} where hypergeometric methods were used instead), can
be found in the appendix together with a certain form of the Lagrange
Inversion Theorem, that is used to give our results a nice form.

\section{Notation}

\subsection{Completed graded super-rings}

By a \emph{(unital) super-ring $R$} we understand a $\set
Z/(2)$-graded unital ring that is commutative in the graded sense. In
particular $R$ comes with a decomposition $R = R_0 \oplus R_1$ such
that in particular $R_0$ is a subring and $R_1$ is a module over
$R_0$. The elements in $R_0$ are the \emph{even elements} (or
\emph{elements of even parity}) and the elements in $R_1$ are the
\emph{odd} ones.

Let $R$ be such a super-ring. A \emph{grading on $R$} is a
decomposition $R = \bigoplus_{n \ge 0} R^n$ of $R$ into subgroups such
that $R^n$ for $n$ even consists entirely of even elements and $R^n$
for $n$ odd consists entirely of odd elements. We require further that
$R^m \cdot R^n \subset R^{m + n}$. A super-ring with a grading is a
\emph{graded super-ring}. The elements in $R^n$ are the
\emph{homogeneous elements of degree $n$}.

A \emph{completed graded super-ring $\hat R$} is a super-ring $\hat R$
that comes with an inclusion $R \subset \hat R$ of a graded super-ring
$R$ such that $\hat R$ is the $I$-adic completion of $R$ with respect to the
ideal $I := \bigoplus_{n > 0} R^n$.

An ordinary ring $R$ becomes trivially a completed graded super-ring
by choosing the grading $R = R^0$. For example, we may view $\set Q$
as a completed graded super-ring. Every graded super-ring can be made
into a completed graded super-ring by passing to its completion with
respect to the ideal generated by the elements of positive degree.

All of the rings that occur in this article will be completed graded
super-ring. In order to simplify notions, we will sloppily call these
completed graded super-rings just \emph{rings}. It frequently happens
that a construction, e.g.~the tensor product, does not yield a
completed graded super-ring but just a graded super-ring. In these
cases, it is to be understood that one has to implicitely pass to the
completion. In other words, ``tensor product'' will always mean
``completed tensor product'', etc.

For the rest of this paper, $R$ will be a ring, i.e.~a completed graded
super-ring.

Let us mention that there are the obvious notions of modules and algebras
over completed graded rings. We will make heavy use of them. One
example will be the cohomology ring $H^*(X, R)$ of a complex manifold
$X$. It is an $R$-algebra (in the above sense) with the grading
$\bigoplus_{p + q \ge 0} H^p(X, R^q)$.

Given an $R$-module $M$ and a natural number $l \in \set N_0$, we
denote by $M[-l]$ the $R$-module whose grading is given by
$$(M[-l])^{n} = \begin{cases}
  M^{n - l} & \text{for $n \ge l$}
  \\
  0 & \text{otherwise}.
\end{cases}$$

\subsection{Partitions}
By $\set P$, we denote the set of partitions of the positive natural
numbers, i.e.~non-increasing finite sequences $(\lambda_1, \dots,
\lambda_n)$ of positive natural numbers. We view the non-zero natural
numbers $\set N$ as a subset of $\set P$ by considering a number $n
\in \set N$ as the trivial partition $(n)$ of itself.

Given a partition $\lambda = (\lambda_1, \dots, \lambda_r) \in \set
P$, we call
$$|\lambda| := r$$
the \emph{length of $\lambda$} and
$$\|\lambda\| := \lambda_1 + \dots + \lambda_r$$
the \emph{weight of
  $\lambda$}, i.e.~a partition $\lambda$ of weight $n$ is a partition
of $n$.

\subsection{The symmetric algebra}

Let $M$ be an $R$-module. We denote by $S_R M$ its (completed)
symmetric algebra (symmetric in the super-sense, of course), i.e.~$S_R
M$ is the free $R$-algebra (in our category of completed graded
commutative unital super-$R$-algebras) generated by $M$. It carries a
natural decomposition
$$S_R M = \prod_{k = 0}^\infty S_R^k M$$
with $S_R^0 M = R$ and
$S_R^1 M = M$.  Given an element $q \in S_R M$, we can write
$$q = \sum_{k = 0}^\infty q_k$$ with
each $q_k \in S_R^k M$. 
For each $n \in R$, we set
$$q(n) := \sum_{k = 0}^\infty n^k q_k.$$

\subsection{Complex manifolds}

By $\emptyset$ we denote the complex manifold without any point. By
$*$ we denote the zero-dimensional complex manifold with exactly one
point.

For each complex manifold $X$, let $\Theta_X$ be its tangent sheaf.

Given any coherent sheaf $\sheaf F$ on $X$, its Chern character in
$H^{2*}(X, \set Q) \subset H^*(X, R)$ is denoted by $\ch(\sheaf F)$.
It decomposes as
$$\ch(\sheaf F) = \sum_{n = 0}^\infty \ch_n(\sheaf F)$$
with $\ch_n(\sheaf F)$ being of degree $2n$.

Given a product $X \times Y$ of two complex manifolds $X$ and $Y$, we
usually identify the cohomology space $H^*(X \times Y, R)$ with the
tensor product $H^*(X, R) \otimes_R H^*(Y, R)$ given by K{\"u}nneth's
theorem.

\section{Characteristic classes}

\subsection{Marked complex manifolds}

\begin{definition}
  A \emph{marking $\alpha$ (over $R$) on a complex manifold $X$} is an
  element $\alpha \in H^*(X, R)$.
  
  A \emph{marked complex manifold $(X, \alpha)$} is a complex
  manifold $X$ together with a marking $\alpha$ on $X$.
  
  A \emph{morphism $f: (X, \alpha) \to (Y, \beta)$ between marked
    complex manifolds $(X, \alpha)$ and $(Y, \beta)$} is a holomorphic
  map $f: X \to Y$ such that $f^*(\beta) = \alpha$.
\end{definition}

Given a marking $\alpha$ on a complex manifold $X$, we write
$$\alpha = \sum_{n = 0}^\infty \alpha_n$$
with $\alpha_n \in
\bigoplus_{p + q = n} H^p(X, R^q)$ being the component of $\alpha$ of
(total) degree $n$.

\begin{definition}
  Let $(X, \alpha)$ and $(Y, \beta)$ be two marked complex
  manifolds.
  
  Their \emph{sum $(X, \alpha) + (Y, \beta)$} is given by $$(X \amalg
  Y, i_* \alpha + j_* \beta),$$
  where $i: X \to X \amalg Y$ and $j: Y
  \to X \amalg Y$ are the natural inclusion maps.
  
  Their \emph{product $(X, \alpha) \cdot (Y, \beta)$} is given by $$(X
  \times Y, p^*\alpha + q^*\beta),$$
  where $p: X \times Y \to X$ and
  $q: X \times Y \to Y$ are the natural projection maps.
\end{definition}

\begin{remark}
  The neutral element with respect to the sum is the marked complex
  manifold $(\emptyset, 0)$. The neutral element with respect to the
  product is the marked complex manifold $(*, 0)$. We thus may
  speak of the \emph{commutative semiring of marked complex manifolds}.
\end{remark}

\begin{remark}
  By formally adding additive inverses, we may turn this semiring into
  the ``ring'' of (virtual) marked complex manifolds. The word
  ``ring'' is in quotation marks as this ``ring'' is a proper class
  and not a set. By tensoring this ring with any other ring $S$, we
  arrive at the $S$-algebra of marked complex manifolds, i.e.~we can study
  formal $S$-linear combinations of marked complex manifolds. All the
  constructions for marked complex manifolds that will follow also extend
  to this formal combinations.
\end{remark}

\subsection{Characteristic classes}

We define the ring
$$U_R := R[p_0, p_1, p_2, \dots, b_0, b_{\frac 1 2}, b_1, \dots]$$
with the grading defined by $\deg p_i = 2 i$ and $\deg b_i = 2 i$.
(Note that this ring is in fact a power series ring as we have to
implicitely complete it.)

\begin{definition}
  The $R$-algebra $U_R$ is the \emph{ring of characteristic
    classes of marked complex manifolds over $R$}.
\end{definition}

The reason for giving this name to the ring is the following: Given a
marked complex manifold $(X, \alpha)$, there is a unique homomorphism
$$U_R \to H^*(X, R),\quad \gamma \mapsto \gamma|_{(X,
  \alpha)}$$
of $R$-algebras with
$$p_i|_{(X, \alpha)} = \ch_i(\Theta_X)\quad\text{and}\quad
b_i|_{(X, \alpha)} = \alpha_{2 i}.$$

We set
$$\ch := \sum_{i = 0}^\infty p_i\quad\text{and}\quad
b := \sum_{j = 0}^\infty b_{j/2}.$$
Then $\ch|_{(X, \alpha)} = \ch(\Theta_X)$ and $b|_{(X, \alpha)} =
\alpha$.
\begin{definition}
  The element $\ch \in U_R$ is the \emph{universal Chern character},
  the element $b \in U_R$ is the \emph{universal marking}.
\end{definition}

\begin{lemma}
  Let $\gamma, \gamma' \in U_R$ with $\gamma|_{(X, \alpha)} =
  \gamma'|_{(X, \alpha)}$ for all marked complex manifolds $(X,
  \alpha)$. Then $\gamma = \gamma'$.
\end{lemma}

\begin{proof}
  This result is well-known for complex manifolds without markings,
  you may consult~\cite{HIR}.

  It is enough to consider the case that $\gamma$ and $\gamma'$ are
  sums of homogeneous elements of degree $n \in \set N$ or less. The
  general case then follows by passing to the limit.  Let us also
  assume for a moment that $\gamma$ and $\gamma'$ do not contain $p_0$
  and $b_0$.

  Consider the (compact) complex manifold
  $$X := \Grass_k(\set C^{n + k}) \times \set C^r/\set Z^{2r},$$
  which
  is the product of the Gra\ss mannian of $k$-dimensional subspaces in
  $\set C^{n + k}$ and an $r$-dimensional torus with $r, k \gg 0$.
  Recall that the rational cohomology of the Gra\ss mannian is freely
  generated up to degree $n$ by the components $\ch_1, \ch_2, \dots,
  \ch_k$ of the Chern character of its tangent bundle. Further, we
  have $H^p(\set C^r/\set Z^{2r}, R) = \Lambda^p_R (\set C^r)$, so we
  find cohomology classes $\alpha_1, \alpha_2, \dots, \alpha_n$ in the
  cohomology of the torus that are algebraic independent up to degree
  $n$, provided $r$ is big enough. These classes define a marking
  $\alpha$ on $X$.

  By these freeness results, we find that $\gamma|_{(X, \alpha)} =
  \gamma'|_{(X, \alpha)}$ implies $\gamma = \gamma'$.
  
  It remains to consider the more general case with $p_0$ and $b_0$
  appearing in $\gamma$. For this, consider $\gamma$ as a polynomial
  $\gamma(p_0, b_0)$ in $p_0$ and $b_0$. Let $(X, \alpha)$ be any
  marked complex manifold.  From this, we construct for each pair $r,
  s \in \set Z$ a new marked complex manifold $(X \times \set
  C^r/\set Z^{2r}, \alpha + s)$ with $\alpha + s$ being the image of
  $\alpha + s \in H^*(X, R)$ under the inclusion map
  $$H^*(X, R) \to H^*(X \times \set C^r/\set Z^{2r}, R) \simeq H^*(X,
  R) \otimes_R H^*(\set C^r/\set Z^{2r}, R), \quad v \mapsto v \otimes
  1.$$
  By virtue of this embedding, we consider $H^*(X, R)$ as a
  subspace of the cohomology ring $H^*(X \times \set C^r/\set Z^{2r},
  R)$. Then we can formulate:
  $$\gamma|_{(X \times \set C^r/\set Z^{2r}, \alpha + s)} = \gamma(p_0
  + r, b_0 + s)|_{(X, \alpha)}.$$
  The same equation holds true for $\gamma'$ in place of
  $\gamma$. Thus, we have
  $$\gamma(p_0 + r, b_0 + s)|_{(X, \alpha)} = \gamma'(p_0 + r, b_0 +
  s)|_{(X, \alpha)}.$$
  Expand both sides in powers of $p_0 + r$ and $b_0 + s$, e.g.~
  $$\gamma(p_0 + r, b_0 + s)|_{(X, \alpha)} = \sum_{l, m \ge 0}
  \gamma_{l, m} (p_0 + r)^l (b_0 + s)^m|_{(X, \alpha)}$$
  with certain
  coefficients $\gamma_{l, m} \in U_R$ that do not contain $p_0$ and
  $b_0$. We get a another expansion with $\gamma$ replaced by
  $\gamma'$. By comparing coefficients, it follows that
  $$\gamma_{l, m}|_{(X, \alpha)} = \gamma'_{l, m}|_{(X, \alpha)}$$ for
  all $l, m \ge 0$. By the previous result on the case where $p_0$ and
  $b_0$ do not appear, we get $\gamma_{l, m} = \gamma'_{l, m}$ and
  therefore $\gamma = \gamma'$.
\end{proof}

\begin{remark}
  Each characteristic class $\gamma \in U_R$ is
  ``additive'' in the following sense: It is
  $$
  \gamma|_{(\emptyset, 0)} = 0
  $$
  and
  $$
  \gamma|_{(X, \alpha) + (Y, \beta)} = i_* \gamma|_{(X, \alpha)} +
  j_* \gamma|_{(Y, \beta)}
  $$
  in $H^*(X \amalg Y, R)$ for all marked manifolds $(X, \alpha)$
  and $(Y, \beta)$.
\end{remark}

\subsection{Multiplicative classes}

On the other hand, not every class is multiplicative in the following
sense:
\begin{definition}
  A characteristic class $\gamma \in U_R$ is
  \emph{multiplicative} if
  \[
  \gamma|_{(\{*\}, 0)} = 1
  \]
  and
  \[
  \gamma|_{(X, \alpha) \cdot (Y, \beta)} = p^* \gamma|_{(X, \alpha)} \otimes
  q^*\gamma|_{(Y, \beta)}
  \]
  in $H^*(X \times Y, R) \simeq H^*(X, R) \otimes_R H^*(Y, R)$ for all
  marked manifolds $(X, \alpha)$ and $(Y, \beta)$. Here as before, $p$
  and $q$ denote the two projections from $X \times Y$ onto its two
  factors.
\end{definition}

To describe all multiplicative classes, we have to consider the
ring $$R[\lambda_1, \lambda_2, \dots]$$
with $\deg \lambda_i = 2$. The symmetric group $$\SG = \colim_{n \in
  \set N_0} \SG_n$$ acts naturally on this $R$-algebra by permuting the
variables $\lambda_i$. We define
$$
\Lambda_R := R[\lambda_1, \lambda_2, \dots]^{\SG},
$$
i.e.~$\Lambda_R$ is the $R$-algebra of $\SG$-invariants of the
polynomial ring.

\begin{definition}
  The $R$-algebra $\Lambda_R$ is the \emph{(completed)
    ring of symmetric functions over $R$.}
\end{definition}

\begin{example}
  The $n$-th elementary symmetric function in the $\lambda_i$ is an
  element in $\Lambda_R$. It is denoted by $c_n$, i.e.~
  $$
  \sum_{i = 0}^\infty c_i = \prod_{i = 1}^\infty (1 + \lambda_i)
  $$
  in $\Lambda_R$ with $\deg c_i = 2i$.
\end{example}

\begin{remark}
  By the theorem of elementary symmetric functions, the homomorphism
  $$
  R[c_1, c_2, \dots] \to \Lambda_R
  $$
  of $R$-algebras that map the (formal) variable $c_i$ of degree
  $2i$ to the $n$-th symmetric function $c_i$ is an isomorphism.
\end{remark}

Consider the (completed) tensor product $\Lambda_R \otimes_R
\Lambda_R$. To distinguish between the first and the second factor,
let us decorate the variables coming from the second factor with a
prime, e.g.~we write $\lambda_i'$. Thus we have
$$\Lambda_R \otimes_R \Lambda_R = R[\lambda_1, \lambda_2, \dots,
\lambda_1', \lambda_2', \dots]^{\SG \times \SG}.$$
Similarly,
everything else arising from the second factor, e.g.~any elementary
symmetric function, is primed.

\begin{remark}
  There is a unique homomorphism
  $$\Lambda_R \otimes_R \Lambda_R \to U_R$$
  of $R$-algebras mapping
  $\sum_{i = 1}^\infty \lambda_i^n$ to $p_n$ and $\sum_{i = 1}^\infty
  \lambda_i'^n$ to $b_n$ for all $n \in \set N$. This homomorphism is
  a monomorphism. By virtue of this homomorphism, we view $\Lambda_R
  \otimes_R \Lambda_R$ as a subring of $U_R$.
\end{remark}

We make use of this identification in the following theorem:  
\begin{theorem}
  Let $f, f' \in 1 + x R[[x]]$ be two power series in $R$ with $f_0 =
  f'_0 = 1$ Then
  \[
  \prod_{i = 1}^\infty f(\lambda_i) \cdot \prod_{i = 1}^\infty
  f'(\lambda_i')
  \]
  is a multiplicative characteristic class in $U_R$.
  
  On the other hand, every multiplicative class in uniquely given by
  such a construction.
\end{theorem}

\begin{proof}
  This result is classical for complex manifolds without marking. See again~\cite{HIR}.

  To see that $\gamma = \prod_{i = 1}^\infty f(\lambda_i) \cdot
  \prod_{i = 1}^\infty f'(\lambda_i')$ is multiplicative, we write
  $$\gamma = \exp\left(\sum_{n \ge 1}(a_n p_n + a_n' b_n)\right)$$
  with $\log f = \sum_{n \ge 1} a_n x^n$ and $\log f' = \sum_{n \ge 1}
  a_n' x^n$. 
  Then one calculates
  $\gamma|_{(*, 0)} = \exp(0) = 1$
  and
  \begin{multline*}
    \gamma|_{(X, \alpha) + (X', \alpha')} = \exp
    \left(\sum_{n \ge 1} \left((a_n p_n + a_n' b_n)|_{(X, \alpha)}
      + (a_n p_n + a_n' b_n)|_{(X', \alpha')}\right)\right)
    \\
    = \gamma|_{(X, \alpha)} \cdot \gamma|_{(X', \alpha')}.
  \end{multline*}

  It remains to prove the implication in the other direction. So let
  $\gamma \in U_R$ be a multiplicative class. From $\gamma|_{(*, 0)} =
  1$ it follows, that the logarithm $\log \gamma \in U_R$ exists. Let
  $$m: H^*(X^n, R) \simeq H^*(X, R)^{\otimes n} \to H^*(X, R)$$
  be the
  multiplication map given the cup-product. As $\gamma$ is
  multiplicative, we have
  $$m\left(\log \gamma|_{(X, \alpha)^n}\right) = n
  \left(\log\gamma|_{(X, \alpha)}\right).$$
  Any power series (with trivial
  constant term) in the $p_i$ and $b_i$ that fulfills this equality in
  place of $\log\gamma$, has to be linear in the $p_i$ and $b_i$. Thus
  $\log \gamma$ is an $R$-linear combination of the $p_i$ and $b_i$,
  which proves the other direction.
\end{proof}

\begin{example}
  The \emph{total Chern classes} $c := 1 + c_1 + c_2 + \dots$ and $c'
  := 1 + c_1' + c_2' + \dots$ are multiplicative classes.
\end{example}

\section{The Hilbert schemes of points on surfaces}

\subsection{Nakajima's creation operators}

Let $X$ be a complex surface. For each $n \in \set N_0$ we denote by
$X^{[n]}$ the Douady space (see~\cite{CM}) of the zero-dimensional subspaces of $X$ of
length $n$, i.e.~holomorphic maps from a complex space $S$ to
$X^{[n]}$ parametrise subspaces $Z$ of $S \times X$ such that the
restriction $p|_Z: Z \to X^{[n]}$ of the projection map $S \times X
\to S$ is flat and finite of degree $n$.

In case $X$ is an algebraic surface, the Douady spaces are algebraic
as well. In fact, they are the corresponding Hilbert schemes of points
on $X$ viewed as a scheme. By abuse of notion, we shall not
distinguish between the two notions:
\begin{definition}
  The space $X^{[n]}$ is the \emph{Hilbert scheme of $n$ points on $X$}.
\end{definition}

\begin{example}
  It is $\Hilb 0 X = *$, the complex manifold consisting of one point,
  and $\Hilb 1 X = X$.
\end{example}

For each point $\xi \in \Hilb n X$, we denote by $\supp \xi = \supp
\sheaf O_\xi$ its \emph{support} which is a finite subset of $X$ whose
cardinality lies between $1$ and $n$. Given two zero-dimensional
subspaces $\xi$ and $\xi'$ of $X$ with $\supp \xi \subset \supp \xi'$,
we denote by $\sheaf I_{\xi, \xi'}$ the kernel of the natural
epimorphism $\sheaf O_{\xi'} \to \sheaf O_{\xi}$.

The identity $X^{[n]} \to X^{[n]}$ corresponds to a subspace $\Xi^n$
of $\Hilb n X \times X$ such that the restriction $p|_{\Xi^n}: \Xi^n
\to \Hilb n X$ of the projection $p: \Hilb n X \times X \to \Hilb n X$
is flat and finite of degree $n$. We denote its structure sheaf by
$\sheaf O_n$ and use the same symbol to denote the corresponding
quotient sheaf of $\sheaf O_{\Hilb n X \times X}$ on $\Hilb n X \times
X$.

\begin{definition}
  The space $\Xi^n$ is the \emph{universal family over $\Hilb n X$}.
\end{definition}

We generalise this subspace as follows: Let $n, l \in \set N_0$ be two
natural numbers. By $\Xi^{n, l}$ we denote the reduced subspace of $\Hilb {l +
  n} X \times \Hilb l X \times X$ whose support is the closure of the
subset
$$\{(\xi', \xi, x) \in \Hilb {l + n} X \times \Hilb l X \times X |
\supp \sheaf I_{\xi, \xi'} \subset \{x\}\}$$

We denote the projections of $\Hilb {l + n} X \times \Hilb l X \times
X$ uniformly by $p$, $q$ and $r$ as given in the diagram
$$
\begin{CD}
  \Hilb {l + n} X @<p<< \Hilb {l + n} X \times \Hilb l X \times X
  @>r>> X \\
  & & @VqVV \\
  & & \Hilb l X.
\end{CD}
$$

The restriction $p|_{\Xi^{n, l}}: \Xi^{n, l} \to \Hilb {n + l} X$ of
$p$ is proper. Thus a push-forward operator
$$p_*(\cdot \cap [\Xi^{n, l}]): H^*(\Hilb {l + n} X \times \Hilb l X
\times X, R) \to H^*(\Hilb {l + n} X, R)$$
is well-defined. We use this map to construct Nakajima's creation
operators in the sequel.

In order to do this, it is convenient to consider all Hilbert schemes of
points of $X$ at once. For this, we define
$$\Hilbtot X := \coprod_{l = 0}^\infty \Hilb l X.$$
The total
cohomology space (up to some shifts to make it completed graded) is
given by
$$\set H_R X := \prod_{l = 0}^\infty H^*(\Hilb l X, R)[-2l].$$
There
are (at least) two distinguished elements in this cohomology space:
One of them is the unit $\one$ of the (unshifted) cohomology ring,
which is given by the infinite sum of the units in each $H^*(\Hilb n
X, R)$, and the other one is the so-called \emph{vacuum $\vac$} that
is the image of $1 \in H^*(\Hilb 0 X, R) = H^*(*, R) = R$ in $\set H_R
X$.

Elements in the subspace $H^*(\Hilb l X, R)[-2l]$ of $H^*(\Hilbtot X,
R)$, $l \in \set N_0$, are called elements \emph{of weight $l$}.

Given $n \in \set N$, let
$$q_n: H^*(X, R) \to \End \set H_R X$$
be the operator-valued linear
map such that for each $\alpha \in H^*(X, R)$ the restriction of
$q_n(\alpha) \in \End \set H_R X$ to $H^*(\Hilb l X, R)$ is given by
$$q_n(\alpha): H^*(\Hilb l X, R) \to H^*(\Hilb {l + n} X, R),\quad \beta
\mapsto p_*((r^* \alpha \cup q^* \beta) \cap [\Xi^{n, l}]).$$
Given $n_1, \dots, n_r \in \set N$, we set
\begin{align*}
  q_{n_1} \circ \dots \circ q_{n_r}: H^*(X, R)^{\otimes r} & \to
  \End \set H_R X,
  \\
  \alpha_1 \otimes \dots \otimes \alpha_{n_r} & \mapsto
  q_{n_1}(\alpha_1) \circ \dots \circ q_{n_r}(\alpha_r),
  \\
  \intertext{and for a partition $\lambda = (\lambda_1, \dots,
    \lambda_r)$ of length $r$, we set} q_\lambda: H^*(X, R) & \to \End
  \set H_R X,
  \\
  \alpha & \mapsto (q_{n_1} \circ \dots \circ q_{n_r})(\delta_*
  \alpha),
\end{align*}
where $\delta: X \to X^r$ is the diagonal map, inducing a push-forward
map $$\delta_*: H^*(X, R) \to H^*(X^r, R) \simeq H^*(X, R)^{\otimes r}$$
on
cohomology.

\begin{definition}
  The operators $q_n(\alpha)$ and $q_\lambda(\alpha)$ for $\alpha \in
  H^*(X, R)$, $n \in \set N$ and $\lambda \in \set P$ are
  \emph{Nakajima's creation operators over $X$ (and $R$)}.
\end{definition}

\begin{remark}
  For homogeneous $\alpha$, the operator $q_\lambda(\alpha)$ is an
  operator of degree $$\deg \alpha + 4 \norm \lambda + 2 \abs \lambda -
  4$$ (and weight $2 \norm \lambda$).
\end{remark}

Set
$$\set N H^*(X, R) := \prod_{n = 1}^\infty H^*(X, R)[-2n].$$
The image
of $\alpha \in H^*(X, R)$ in $\set N H^*(X, R)$ under the inclusion
map corresponding to a positive natural number $n$ is denoted by
$q_n(\alpha)$, which is not to be confused with the \emph{operator}
$q_n(\alpha) \in \End \set H_R X$. Finally consider the symmetric algebra
$$N_R X := S_R (\set N H^*(X, R)).$$ (The ``N'' stands for Nakajima.)
\begin{theorem}
  The $R$-linear map
  $$N_R X \to \set H_R X,\quad q_{n_1}(\alpha_1) \dots
  q_{n_r}(\alpha_r) \mapsto (q_{n_1}(\alpha_1) \circ \dots \circ
  q_{n_r}(\alpha_r)) \vac
  $$
  is an isomorphism of $R$-modules.
  
  In particular, any two creation operators commute (in the
  super-sense, of course).
\end{theorem}

\begin{proof}
  The reader may find the original proof in~\cite{NAK}.
\end{proof}

\begin{example}
  It is $$\one = \exp(q_1(1)) \vac.$$
\end{example}

The multiplication morphism $(N_R X)^{\otimes n} \to N_R X$ of the
$R$-algebra $N_R X$ induces via the isomorphism in Nakajima's theorem
a ``multiplication map''
$$\mu: (\set H_R X)^{\otimes n} \to \set H_R X,$$
which is \emph{not} the cup-product.

We can express the Poincar\'e pairing on the cohomology of the Hilbert
schemes of points on a compact surface $X$ in terms of the creation
operators as follows:
\begin{proposition}
  Given $\alpha_1, \dots, \alpha_r, \beta_1, \dots, \beta_s \in H^*(X,
  R)$ and natural numbers $m_1, \dots, m_r, n_1, \dots, n_s \in \set N$ we do have
  \begin{multline*}
    \langle q_{m_1}(\alpha_1) \dots q_{m_r}(\alpha_r)\vac,
    q_{n_1}(\beta_1) \dots q_{n_s}(\beta_r)\vac \rangle
    \\
    = \sum_{\substack{\sigma: \{1, \dots, s\}\\ \tilde \to \{1, \dots,
        r\}}} \pm \prod_{i = 1}^s m_i \cdot \delta_{m_i,
      n_{\sigma(i)}} \cdot \langle\alpha_i, \beta_{\sigma(i)}\rangle
  \end{multline*}
  where the brackets denote the Poincar\'e pairing on $\set H_R X$ and
  $H^*(X, R)$, respectively, and the sum runs over all
  \emph{bijections}. The sign is the Koszul sign arising from
  reordering the $\alpha_i$ and $\beta_i$.
  
  In particular, the pairing vanishes if $r \neq s$ or each index of a
  creation operator does not occur with the same multiplicity on both
  sides.
\end{proposition}

\begin{proof}
  There are adjoint operators to the creation operators acting on
  $\set H_R X$, called \emph{annihilation operators}. Nakajima
  calculated commutators, which are fundamental for the study of the
  cohomology of the Hilbert schemes, between the creation and
  annihilation operators. From these commutation relations, one can
  deduce the above result. The way to do this can be found
  in~\cite{LHS} (see in particular the first paragraph in section 6).
\end{proof}

In the sequel we need a simple fact following from this description of
the Poincar\'e pairing. This fact is stated in the following lemma
where we also assume that $X$ is compact. Let further $\alpha_1,
\dots, \alpha_n \in H^*(X, R)$ be linearly independent cohomology
classes on the surface. We choose classes $\beta_1, \dots, \beta_n$
such that $\langle\alpha_i, \beta_j\rangle = \delta_{ij}$ for $i, j
\in \{1, \dots, n\}$. For each $r \ge 1$ we set $\check \beta_i^r :=
p^* \beta$ where $p: X^r \to X$ is the projection onto the first
factor.  Recall the push-forward map $\delta_*: H^*(X, R) \to H^*(X^r,
R)$. We set $\alpha_i^r := \delta_*(\alpha_i)$. One verifies that
$\langle\check\beta_i^s, \alpha_j^r \rangle \neq 0$ if $i = j$ and $s
= r$ and zero otherwise. For each partition $\lambda = (n_1, \dots,
n_r)$, we finally set $$\check q_\lambda(\beta_i) := (q_{n_1} \circ
\dots \circ q_{n_r})(\check\beta_i^r)$$
for $i \in \{1, \dots, n\}$.
\begin{lemma}
  \label{lem:orthogonality2}
  Given partitions $\lambda_1, \dots, \lambda_r, \mu_1, \dots, \mu_s$
  and maps $\rho: \{1, \dots, r\} \to \{1, \dots, n\}$ and $\sigma:
  \{1, \dots, s\} \to \{1, \dots, n\}$, the expression
  $$\langle \check q_{\lambda_1}(\beta_{\rho(1)}) \dots \check
  q_{\lambda_r}(\beta_{\rho(r)}) \vac,
  q_{\mu_1}(\alpha_{\sigma(1)})\dots
  q_{\mu_r}(\alpha_{\sigma(r)})\vac\rangle$$
  is non-zero if and only
  there is a bijection $\tau: \{1, \dots, r\} \to \{1, \dots, s\}$
  such that $\lambda_i = \mu_{\tau(i)}$ and $\rho(i) =
  \sigma(\tau(i))$.
\end{lemma}

\begin{proof}
  This follows from the previous proposition and the orthogonality of
  the $\alpha_i$ and $\beta_i$.
\end{proof}

\subsection{Hilbert schemes of marked complex surfaces}

For each marking $\alpha \in H^*(X, R)$ on a complex surface $X$ and
each $n \in \set N_0$, we set
\[
\alpha^{[n]} := p_*(\ch(\sheaf O_{\Xi^n}) \cup q^*(\td(X) \cup \alpha)),
\]
which is a marking on $\Hilb n X$, which depends $R$-linearly on
$\alpha$.  The corresponding marked complex manifold is denoted by
$$(X, \alpha)^{[n]} := (X^{[n]}, \alpha^{[n]}).$$

\begin{definition}
  The marked complex manifold $(X, \alpha)^{[n]}$ is \emph{the Hilbert
    scheme of $n$ points on $(X, \alpha)$}.
\end{definition}

Again, it is a good idea to consider all these spaces at once, so we set
\[
\Hilbtot (X, \alpha) := \sum_{l = 0}^\infty (X, \alpha)^{[n]}.
\]

The function $\Hilbtot(\cdot)$ works like an exponential. By this, we
mean the following:
\begin{proposition}
  Let $(X, \alpha)$ and $(Y, \beta)$ be two marked complex
  surfaces. Then
  \begin{align*}
    \Hilbtot(\emptyset, 0) & = (*, 0)\\
    \intertext{and}
    \Hilbtot((X, \alpha) + (Y, \beta)) & = \Hilbtot(X, \alpha)
    \cdot \Hilbtot (Y, \beta).
  \end{align*}
\end{proposition}

\begin{proof}
  From the fact that a zero-dimensional subspace of $X \amalg Y$ of
  length $n$ is given by a pair of a zero-dimensional subspace of $X$
  of length $p$ and a zero-dimensional subspace of $Y$ of length $q$
  with $p + q = n$, one proves that
  $$\Hilb n {(X \amalg Y)} = \coprod_{p + q = n} \Hilb p X \times
  \Hilb q Y,$$
  see e.g.~\cite{EGL}.
  Furthermore, the universal family over $\Hilb n {(X \amalg Y)}$
  restricted to the component $\Hilb p X \times \Hilb q Y$ is given by
  $$(\Xi^p \times \Hilb q Y) \amalg (\Hilb p X \times \Upsilon^q)
  \subset (\Hilb p X \times X \times \Hilb q Y) \amalg (\Hilb p X
  \times \Hilb q Y \times Y),$$
  where $\Xi^p$ is the universal family
  over $\Hilb p X$ and $\Upsilon^q$ is the universal family over
  $\Hilb q Y$. Thus by definition of $(\alpha + \beta)^{[n]}$, its
  restriction onto the component $\Hilb p X \times \Hilb q Y$ is given
  by
  $$\alpha^{[p]} \boxtimes 1 + 1 \boxtimes \beta^{[q]} \in H^*(\Hilb p
  X \times \Hilb q Y, R).$$
  Summing all up, this yields the claimed formula.
\end{proof}

Let $n \in \set N_0$. Given a marked complex surface $(X, \alpha)$, we
may consider the $n$-fold sum $n (X, \alpha)$, which is again a marked
surface. The previous proposition tells us that
$$\Hilbtot(n (X, \alpha)) = \Hilbtot(X, \alpha)^n.$$
In particular, we have
$$\set H_R (n (X, \alpha)) = (\set H_R (X, \alpha))^{\otimes n}.$$
Recall the $R$-linear map $\mu: (\set H_R X)^{\otimes n} \to \set H_R
X$ given by Nakajima's description of the cohomology space in terms of
the creation operators applied to the vacuum. This give us thus an
$R$-linear map $$\mu: \set H_R (n (X, \alpha)) \to \set H_R (X,
\alpha).$$
We shall need this map in the next section.

\section{Universal formulas}

Consider the $\set Q$-algebra
$$
A := \set Q[K, e, \alpha_0, \alpha_1, \alpha_2, \alpha_3, \alpha_4]
$$
with $\deg K = 2$, $\deg e = 4$ and $\deg \alpha_i = i$. Let $A^{\ge 5}$
be the ideal generated by the elements of degree five or more in $A$.
We set
$$
V_R := R \otimes_{\set Q} (A/A^{\ge 5}).
$$
The element $K \in V_R$ is called the \emph{universal canonical class}, the
element $e \in V_R$ is called the \emph{universal Euler class}.

\begin{definition}
  The $R$-algebra $V_R$ is the \emph{universal cohomology ring
    for marked complex surfaces over $R$}.
\end{definition}

This is due to the following: Given a marked complex surface $(X,
\alpha)$, there is a unique homomorphism
$$V_R \to H^*(X, R),\quad \gamma \mapsto \gamma|_{(X, \alpha)}$$
of
$R$-algebras with
$$K|_{(X, \alpha)} = - c_1(\Theta_X),\quad
e|_{(X, \alpha)} = c_2(\Theta_X),\quad\text{and}\quad
\alpha_i|_{(X, \alpha)} = \alpha_i.$$

Consider a $\set Q$-algebra $S$. A (virtual) marked compact complex
surface $(X, \alpha)$ in the $S$-algebra of marked complex surfaces is called
\emph{versal} if the homomorphism
$$V_R \to H^*(X, R),\quad \mapsto \gamma|_{(X, \alpha)}$$
is a monomorphism. There are versal marked surfaces. They can be constructed
by taking a sufficiently general linear combination of copies of $\set
P^2$, $\set P^1 \times \set P^1$, $\set C^2/\set Z^4$ and $\set P^1
\times \set C/\set Z^2$ with sufficiently general markings.

\begin{lemma}
  Let $\gamma, \gamma' \in V_R$ with $\gamma|_{(X, \alpha)} =
  \gamma'|_{(X, \alpha)}$ for all marked complex surfaces $(X,
  \alpha)$. Then $\gamma = \gamma'$.
\end{lemma}

\begin{proof}
  This follows by evaluating $\gamma$ and $\gamma'$ on a versal marked
  surface.
\end{proof}

Consider the $R$-algebra
$$\set P V_R := \prod_{\lambda \in \set P} V_R[-2\norm \lambda].$$
The canonical injection from $V \to \set P V_R$
corresponding to the factor indexed by a partition $\lambda \in \set P$ is
denoted by
$$V_R \to \set P V_R,\quad v \mapsto q_\lambda(v).$$

An element in $\set P V_R$ is called an element \emph{of weight $l$},
$l \in \set N_0$, if it is a (finite) sum of elements of the form
$q_\lambda(v)$ with $\|\lambda\| = l$. 

We set $$Q_R := S_R(\set P V_R).$$ Note that $\set P V_R$ becomes a
subspace of $Q_R$.

\begin{definition}
  The $R$-algebra $Q_R$ is the
  \emph{universal creation operator algebra of the Hilbert schemes of
    points of marked surfaces over $R$}.
\end{definition}

The algebra $Q_R$ is universal in the following sense: Given a marked
complex surface $(X, \alpha)$, there is a unique homomorphism of
$R$-algebras
$$Q_R \to \End(\set H_R X),\quad q \mapsto q|_{(X, \alpha)}$$
mapping $q_\lambda(v)$ for $v \in V_R$ and $\lambda \in \set P$ to
$q_\lambda(v|_{(X, \alpha)})$.

\begin{lemma}
  \label{lem:uniqueness_in_Q_R}
  Let $q, q' \in Q_R$ with $q|_{(X, \alpha)} = q'|_{(X, \alpha)}$
  for all marked complex surfaces $(X, \alpha)$. Then $q = q'$.
\end{lemma}

\begin{proof}
  We may assume that $X$ is a versal surface. Then $(\alpha_0, \dots,
  \alpha_4, 1, K, e) \in V_R$ maps to a linearly independent system on
  $H^*(X, R)$. We choose a ``dual system'' $(\beta_0, \dots, \beta_7)$
  as in lemma~\ref{lem:orthogonality2} to this system.
  Applying this lemma enables us to extract the components of $q$ and
  $q'$ that have to occur with the same coefficients. From this, we
  conclude $q = q'$.
\end{proof}

Recall the definition of $q(n)$ given in the section about notations.
\begin{lemma}
  Let $(X, \alpha)$ be a marked complex surface and $n \in \set N_0$
  a natural number. Then
  $$
  \mu \circ q|_{n (X, \alpha)} = q(n)|_{(X, \alpha)} \circ \mu
  $$
  as operators from $(\set H_R X)^{\otimes n}$ to $\set H_R X$ for
  all $q \in Q_R$.
\end{lemma}

\begin{proof}
  By linearity, we may assume that $$q
  = q_{\lambda_1}(v_1) \dots q_{\lambda_r}(v_r)$$ for
  partitions $\lambda_1, \dots, \lambda_r \in \set P$ and elements
  $v_1, \dots, v_r \in V_R$.
  Then
  $$\mu(q_{\lambda_1}(v_1) \dots q_{\lambda_r}(v_r)|_{n (X, \alpha)}
  \vac)
  = \mu(q_{\lambda_1}(n \cdot \beta_1) \dots q_{\lambda_r}(n
  \cdot \beta_r)\vac),$$
  with $\beta_i := v_i|_{(X, \alpha)}$ and where
  $$n \cdot \beta_i := \sum_{j = 1}^n (i_j)_* \beta\quad \in H^*(n (X,
  \alpha), R) \simeq H^*(X, \alpha)^{\oplus n}$$
  where $i_j: H^*(X, R)
  \to H^*(X, R)^{\oplus n}$ is the injection map into the $j$-th
  summand. Of course, $\vac$ is the vacuum vector in $\set H_{n (X,
    \alpha)}$.
  
  The operator $\mu$ ``commutes'' with the $q_{\lambda}$ in the sense
  that $$\mu \circ q_{\lambda_i}(n \cdot \beta_i) = n \cdot
  q_{\lambda_i}(\beta_i) \circ \mu.$$
  (This follows from $$q_l(n \cdot
  \beta_i)(\epsilon_1 \otimes \dots \otimes \epsilon_n) = \sum_{i =
    1}^n \epsilon_1 \otimes \dots \otimes \epsilon_{i - 1} \otimes
  q_l(\beta_i)(\epsilon_i) \otimes \epsilon_{i + 1} \otimes \dots
  \otimes \epsilon_n$$
  for $l \in \set N$ and $\epsilon_1 \otimes
  \dots \otimes \epsilon_n \in \set H_R(n \cdot (X, \alpha)) \simeq
  (\set H_R X)^{\otimes n}$.)  

  This gives
  \begin{multline*}
    \mu(q_{\lambda_1}(n \cdot \beta_1) \dots q_{\lambda_r}(n
    \cdot \beta_r)\vac) = n^r (q_{\lambda_1}(\beta_1) \dots
    q_{\lambda_r}(\beta_r))(\mu(\vac))
    \\
    = n^r (q_{\lambda_1}(v_1) \dots
    q_{\lambda_r}(v_r))|_{(X, \alpha)}(\mu(\vac)).
  \end{multline*}
  Together with the fact that acting on $\vac$ is faithful, this yields
  $$\mu \circ q_{\lambda_1}(v_1) \dots q_{\lambda_r}(v_r)|_{n (X,
    \alpha)} = n^r (q_{\lambda_1}(v_1) \dots q_{\lambda_r}(v_r))|_{(X,
    \alpha)} \circ \mu.$$
\end{proof}

\begin{lemma}
  Let $q \in Q_R$ with
  $$q_{(X, \alpha) + (X', \alpha')} = q|_{(X, \alpha)} \otimes \id +
  \id \otimes q|_{(X', \alpha')}$$
  in $\End \set H_R ((X, \alpha) +
  (X', \alpha')) = \End (\set H_R X \otimes \set H_R X')$
  for two
  marked complex surfaces. Then $q \in \set P V_R \subset Q_R$.
\end{lemma}

\begin{proof}
  For each natural number $n \in \set N$ and each complex
  surface $(X, \alpha)$, we have
  $$
  n q|_{(X, \alpha)} \circ \mu = \mu \circ q|_{n (X, \alpha)} =
  q(n)|_{(X, \alpha)} \circ \mu
  $$
  as maps from $(\set H_R X)^{\otimes n} \to \set H_R (X)$, where
  the first equality is due to the assumption and the second due to
  the previous lemma.  It follows that $n q|_{(X, \alpha)} =
  q(n)|_{(X, \alpha)}$ and thus $n q = q(n)$ for all $n \in \set N$.
  This means that $q(n)$ is linear in $n$, i.e.~$q \in \set P V_R$.
\end{proof}

\begin{theorem}
  For each characteristic class $\gamma$ of marked complex
  manifolds, there exists a unique class $q(\gamma) \in Q_R$ such
  that
  \[
  q(\gamma)|_{(X, \alpha)}\one = \gamma|_{\Hilbtot (X, \alpha)}
  \]
  for all marked complex surfaces $(X, \alpha)$.
\end{theorem}

\begin{proof}
  In~\cite{BOS2}, it is stated that each polynomial $\pi$ in the Chern
  classes of tautological sheaves $u^{[\cdot]}$ of the Hilbert
  schemes of points on a surfaces can be expressed by a ``universal
  formula''. In our language, this means that there exists a $q \in
  Q_R$ such that
  $$\sum_{n \ge 0} \pi(u^{[n]}) = q|_{(X, \ch(u))}\vac$$
  for all complex surfaces $X$ and each ``sheaf'' $u \in K(X)$. (In
  fact, the results in~\cite{BOS2} are only stated for algebraic
  surfaces, however, the result generalises easily.)

  We may generalise this result even further: The theorem
  in~\cite{BOS2} remains true if we substitute $\ch(u)$ by an
  arbitrary class $\alpha \in H^*(X, R)$ (however, one has to start to
  keep track of the ``Koszul signs'' arising), i.e.~for each
  polynomial $\pi$ in the cohomology classes $\alpha^{[\cdot]}$,
  there exists a $q \in Q_R$ such that
  $$\sum_{n \ge 0} \pi(\alpha^{[n]}) = q|_{(X, \alpha)}\vac$$
  for each marked complex surface $(X, \alpha)$:
  
  An analogous statement holds true for polynomials in the Chern
  classes of the tangent sheaves of the Hilbert schemes of points on
  a surface, which is also proven in~\cite{BOS2}. Thus, everything
  holds true for ``mixed'' polynomials in the $\alpha^{[\cdot]}$ and
  the Chern classes of the tangent sheaves, and thus, via completion,
  for $\gamma$.

  This proves that there exists a $q(\gamma)$ with
  $$q(\gamma)|_{(X, \alpha)}\one = \gamma_{\Hilbtot (X, \alpha)}.$$
  (Recall that $\one = \exp(q_1(1))\vac$.) The uniqueness of
  $q(\gamma)$ follows from lemma~\ref{lem:uniqueness_in_Q_R}.
\end{proof}

\begin{corollary}
  Consider the universal Chern character $\ch \in U_R$ and the
  universal marking $\alpha \in U_R$. 
  It is $q(\ch) \in \set P V_R \subset Q_R$ and $q(\alpha) \in
  \set P V_R \subset Q_R$.
\end{corollary}

\begin{proof}
  We apply $q(\ch)|_{(X, \alpha) + (X', \alpha')}$ to the unit, which
  gives:
  \begin{align*}
    q(\ch)|_{(X, \alpha) + (X', \alpha')} \one
    & = \ch(\Hilbtot((X, \alpha) + (X', \alpha'))) \\
    & = \ch(\Hilbtot (X, \alpha) \cdot \Hilbtot (X', \alpha')) \\
    & = \ch(\Hilbtot (X, \alpha)) \otimes \one
    + \one \otimes \ch(\Hilbtot (X', \alpha')) \\
    & = \left(q(\ch)|_{(X, \alpha)} \otimes \id
    + \id \otimes q(\ch)|_{(X', \alpha')}\right) \one.
  \end{align*}
  This proves the statement for $q(\ch)$ by the previous lemma.
  
  The proof of the statement for $q(\alpha)$ is virtually the same, so
  we shall omit it.
\end{proof}

\begin{corollary}
  \label{cor:multiplicative_formula}
  For each multiplicative characteristic class $\gamma$ of marked
  complex manifolds, there exists a unique class $\log q(\gamma) \in
  \set P V_R$
  such that
  $$\exp(\log q(\gamma)) = q(\gamma)$$ in $Q_R$.
\end{corollary}

\begin{proof}
  It is
  \[
  q(\gamma)|_{(\emptyset, 0)} \one = \gamma|_{\Hilbtot (\emptyset, 0)}
  = \gamma|_{(*, 0)} = 1,
  \]
  thus $q(\gamma) \in 1 + (\set P V_R) \cdot Q_R \subset Q_R$
  Therefore, there exists a $\log q(\gamma) \in Q_R$ with $\exp(\log
  q(\gamma)) = q(\gamma)$.
  
  It remains to show that $\log q(\gamma) \in \set P V_R \subset Q_R$.
  This follows from the previous lemma since
  $$\log q(\gamma)|_{(X, \alpha) + (X', \alpha')} = \log
  q(\gamma)|_{(X, \alpha)} \otimes \id + \id \otimes \log
  q(\gamma)|_{(X', \alpha')}.$$
  This in turn is due to the fact that
  \begin{align*}
    q(\gamma)|_{(X, \alpha) + (X', \alpha')}\one
    & = \gamma|_{\Hilbtot((X, \alpha) + (X', \alpha'))}\\
    & = \gamma|_{\Hilbtot(X, \alpha) \cdot \Hilbtot(X', \alpha')} \\
    & = \gamma|_{\Hilbtot(X, \alpha)} \otimes \gamma|_{\Hilbtot(X', \alpha')}\\
    & = q(\gamma)|_{(X, \alpha)} \one \otimes q(\gamma)|_{(X',
      \alpha')}
    \one\\
    & = \left(q(\gamma)|_{(X, \alpha)} \otimes q(\gamma)|_{(X',
        \alpha')}\right)\one
  \end{align*}
  for marked complex surfaces $(X, \alpha)$ and $(X', \alpha')$.
\end{proof}

Let $L: V_R \to V_R$ be the homomorphism of $R$-algebras that maps $K$
to $K$, $e$ to $e$, $\alpha_{2i}$ to
$\frac{(\alpha_2)^i}{i!}$ for $i \in \{0, 1, 2\}$ and $\alpha_1$ and
$\alpha_3$ to zero. In particular
$$L(\alpha) = \exp (\alpha_2)$$
with $\alpha = \sum_{i = 0}^4
\alpha_i$.  This homomorphism induces homomorphisms of $\set P V_R$
and $Q_R$, which will also be denoted by $L$.
\begin{proposition}
  Using the homomorphism $L$, we can formulate the following partial
  result on $q(c') \in Q_R$, namely
  $$L(\log q(c')) = \sum_{m = 1}^\infty \frac{(-1)^{m - 1}} m q_m(1
  + \alpha_2) - q_1(1).$$
\end{proposition}

\begin{proof}
  It is
  $$L(\log q(c'))|_{(X, \alpha)}\one = \log q(c')|_{(X, \exp
    \delta)}\one = \sum_{n \ge 0} c'|_{\Hilb n {(X, \alpha)}}$$
  with
  $\delta := \exp(\alpha_2)$. We may assume that $\alpha_2$ is
  algebraic, i.e.~$\alpha_2 = c_1(\sheaf L)$ for a invertible sheaf
  $\sheaf L$ (from which $\delta = \ch(\sheaf L)$ follows). Then the
  proposition follows from Lehn's theorem
  $$\sum_{n \ge 0} c(\sheaf L^{[n]}) = \exp\left(\sum_{m = 1}^\infty
    \frac{(-1)^{m - 1}} m q_m(c_1(\sheaf L))\right) \vac.
  $$
  given in~\cite{LEH}
  as $c(\sheaf L^{[n]}) = c'|(\alpha^{[n]})$.
\end{proof}

\begin{corollary}
  \label{cor:lehn_formula}
  Given a marked complex surface $(X, \alpha)$ with $\alpha =
  \exp(\delta)$ for a class $\delta \in H^2(X, R)$, we have
  $$c'|_{\Hilbtot(X, \alpha)} = \exp\left(\sum_{m = 1}^\infty
    \frac{(-1)^{m - 1}} m q_m(1 + \delta)\right) \vac.
  $$
\end{corollary}
\qed

\section{The Hilbert scheme of the affine plane}

Let $\gamma \in V_R$ be a multiplicative characteristic class,
i.e.~$\gamma = \prod_{i = 1}^\infty f(\lambda_i)$ for a uniquely
determined power series $f \in 1 + x R[[x]]$.

We define a power series $g \in t R[[t]]$ by
$$
\frac {\partial g}{\partial t} \left(\frac x {f(x) f(-x)}\right) =
f(x) f(-x),
$$
i.e.~$t \frac {\partial g(t)} {\partial t}$ is the compositional
inverse of the power series $\frac x {f(x) f(-x)}$. In particular, $g$
is an odd power series.

Let $A: V_R \to V_R$ be the homomorphism of $R$-algebras that maps $K$
and $e$ to zero and $\alpha$ to $1$. The $R$-linear map $\set P V_R
\to \set P V_R$ that maps $q_\lambda(v)$ to $q_\lambda(A(v))$ for
$|\lambda| = 1$ and $q_\lambda(v)$ to zero for $|\lambda| > 1$ will
also be denoted by $A$. This map in turn induces a homomorphism of
$R$-algebra $A: Q_R \to Q_R$.

\begin{theorem}
  Using the homomorphism $A$, the following partial result on
  $q(\gamma) \in Q_R$ holds, namely
  $$
  A(\log q(\gamma)) = \sum_{k \ge 1}^\infty \left(\frac 1 {k!}
    \left.\frac{\partial^k}{\partial t^k}\right|_{t = 0} g
  \right) \cdot
  q_k(1) - q_1(1).$$
\end{theorem}

\begin{proof}
  The proof follows the strategy also used in~\cite{BOS2}. First of
  all note that
  $$A(q(\gamma))|_{X, \alpha}\one = q(\gamma)|_{(\set C^2,
    0)}\one = \sum_{n \ge 0} \gamma|_{\Hilb n{(\set C^2, 0)}}.$$
  Thus the proof of
  the theorem reduces to calculating multiplicative classes on the
  Hilbert schemes of points on the affine plane. By
  corollary~\ref{cor:multiplicative_formula}, the right hand side can
  be written as
  $$\sum_{n \ge 0} \gamma|_{\Hilb n {(\set C^2, 0)}} = 
  \exp\left(\sum_{k \ge 1} g_k q_k(1)\right)\vac$$
  for certain $g_k
  \in R$. The reader may wonder why no $q_{\lambda}(1)$ with
  $|\lambda| \ge 2$ appears on the right hand side, but these
  $q_{\lambda}(1)$ act trivially on $\set H_{\set C^2}$ as $\delta_*:
  H^*(\set C^2, R) \to H^*((\set C^2)^r, R)$ is the zero map for $r
  \ge 2$.
  
  Thus it remains to calculate the $g_k$ that appear above, in fact to
  prove that $g = \sum_{k \ge 1} g_k t^k$. Restricting the previous
  formula to the classes of (maximal) degree $n - 1$ on each Hilbert
  $\Hilb n {(\set C^2)}$ scheme yields
  $$\sum_{k \ge 1} g_k q_k(1)\vac = \sum_{n \ge 0} (\gamma_{n -
    1})|_{\Hilb n{(\set C^2, 0)}}$$
  where $\gamma_{n - 1}$ is the
  component of $\gamma$ of degree $n - 1$. As explained
  in~\cite{BOS2}, this follows from the fact that $H^r(\Hilb n {(\set
    C^2)}, R) = 0$ for $r \ge n$.

  Looking in~\cite{BOS1}, we find the following formula in the proof
  of theorem 5.1:
  \begin{multline*}
    c_{n - 1}|_{\Hilb n {(\set C^2, 0)}} \\
= \sum_{\norm\lambda = n} \frac 1
    {h(\lambda)} \left(\left.\frac 1 {(n - 1)!} \frac{\partial^{n -
            1}}{\partial x^{n - 1}}\right|_{x = 0} \prod_{w \in
        D(\lambda)} ((1 + h(w) x) (1 - h(w) x))\right) \\ \cdot
    \chi^\lambda_{(n)} \cdot z_{(n)}^{-1} \cdot q_n(1)\vac,
  \end{multline*}
  where the symbols $D(\lambda)$, $h(\lambda)$, $h(x)$,
  $\chi_\mu^\lambda$ and $z_\mu$ are defined in~\cite{BOS1}
  (or~\cite{BOS2}). (This formula follows from theorem 4.2
  in~\cite{BOS2}.)

  Slightly generalising this formula by replacing the Chern class
  $c_{n - 1}$ by $\gamma_{n - 1}$, which amounts to replace the power
  series $1 + x$ by the more general $f$, yields
  \begin{multline*}
    \gamma_{n - 1}|_{\Hilb n {(\set C^2, 0)}} = \sum_{\norm\lambda = n} \frac 1
    {h(\lambda)} \left(\left.\frac 1 {(n - 1)!} \frac{\partial^{n -
            1}}{\partial x^{n - 1}}\right|_{x = 0}
      \prod_{w \in
        D(\lambda)} (f(h(w) x) f(-h(w) x))\right)
    \\ \cdot \chi^\lambda_{(n)} \cdot z_{(n)}^{-1} \cdot
    q_n(1)\vac.
  \end{multline*}
  (Already in~\cite{BOS2} is has been stated that similar formulas hold
  true for more than just the total Chern class.)

  We proceed as in~\cite{BOS2}, noting that
  $$
  \chi_{(n)}^\lambda = \begin{cases}
    (-1)^s & \text{for $\lambda = (n - s, \underbrace{1, \dots,
        1}_{\text{$s$-times}})$ and $0 \leq s
      < n$}
    \\
    0 & \text{otherwise}
  \end{cases}
  $$
  and that $z_{(n)} = n$, $\{h(w): w \in
  D(\lambda)\} = \{1, \dots, s, 1, \dots, n - s - 1, n\}$ (with
  multiplicities) and $h(\lambda) = s! \cdot (n - s - 1)! n$ for a partition
  $\lambda$ of $n$ the form $\lambda = (n - s, 1, \dots, 1)$, $0 \leq
  s < n$. This simplifies the above expression for $\gamma_{n -
    1}|_{(\set C^2, 0)}$ to
  \begin{align*}
    \gamma_{n - 1}|_{\Hilb n {(\set C^2, 0)}} & = \frac{q_n(1)} n \vac \cdot \sum_{s =
      0}^{n - 1} \frac{(-1)^s}{n \cdot s!(n - 1 - s)!}
    \\
    & \cdot \left(\left. \frac 1{(n - 1)!}  \frac{\partial^{n -
            1}}{\partial x^{n - 1}}\right|_{x = 0} \prod_{\substack{k
          \in \{1, \dots, s, \\1, \dots, n - 1 - s,
          n\}}} (f(k \cdot x) \cdot f(- k \cdot x))
    \right)
    \\
    & = \frac{q_n(1)} n \vac \cdot
      \left. \frac 1{(n - 1)!} \frac{\partial^{n -
            1}}{\partial x^{n - 1}}\right|_{x = 0} 
      f(n \cdot x) \cdot f(-n \cdot x)
      \\
      & \qquad \cdot
      \sum_{s = 0}^{n - 1} \frac{(-1)^s}{n \cdot s!(n - 1 - s)!}
      \prod_{k = - (n - s)}^s (f(k \cdot x) f(- k \cdot x)).
  \end{align*}
  By lemma~\ref{lem:identity2} (to be found in the appendix)
  this gives
  $$\gamma_{n - 1}|_{\Hilb n {(\set C^2, 0)}} = \frac{q_n(1)} n \vac
  \cdot \left.\frac {(-1)^{n - 1}} {n!}
    \frac{\partial^{n - 1}}{\partial x^{n - 1}}\right|_{x = 0} (f(x) \cdot
  f(-x))^n.$$
  As $f(x) \cdot f(-x)$ is an even power series, we can leave out the
  sign $(-1)^{n - 1}$. Together with the Langrange Inversion Theorem,
  see proposition~\ref{prop:lagrange_expansion} in the appendix, this
  yields
  $$\sum_{n \ge 0} \gamma_{n - 1}|_{\Hilb n {(\set C^2, 0)}}
  = \sum_{k \ge 1} \left(\left.\frac{1}{k!}\frac{\partial^k}{\partial
        t^k}\right|_{t = 0} g\right) \cdot q_k(1) \vac,$$
  which is all that remained to show.
\end{proof}

\begin{corollary}
  Consider the marked complex surface $(\set C^2, 0)$. We have
  $$
  \gamma|_{\Hilbtot (\mathbf C^2, 0)} = \exp\left(
    \sum_{k \ge 1}^\infty \left(\frac 1 {k!}
      \left.\frac{\partial^k}{\partial t^k}\right|_{t = 0} g(t)\right) \cdot
    q_k(1)\right)\vac.
  $$
\end{corollary}
\qed

\begin{example}
  Let $\gamma = c$ be the \emph{total Chern class}. Then $f = 1 + x$
  and thus
  $$g = \frac {1 - \sqrt{1 + 4 t^2}}{2t} + \arcsin(2t) = \sum_{n =
    0}^\infty (-1)^n \frac{1}{(n + 1) \cdot (2n + 1)} \cdot \binom{2n} n \cdot
  t^{2n + 1}.$$
  By the previous theorem, this gives
  $$A(\log q(c)) = \sum_{n \ge 0} \frac{(-1)^n} {n + 1} \cdot
  \binom{2n} n \cdot \frac{q_{2n + 1}(1)}{2n + 1} - q_1(1).$$
  This formula is one of the main results in~\cite{BOS1}.
\end{example}

\begin{example}
  Let $\gamma = s = c^{-1}$ be the \emph{total Segre class}, i.e.~$f = (1
  + x)^{-1}$. Thus
  $$\frac{\partial g}{\partial t}(x \cdot (1 - x^2)) = (1 - x^2)^{-1}.$$
  By the Lagrange Inversion Theorem, we have
  $$g = \sum_{n = 0}^\infty \frac 1 {(2 n + 1)^2} \cdot \binom{3n} n \cdot
  t^{2n + 1},$$  
  and therefore
  $$A(\log q(s)) = \sum_{n \ge 0} \frac{1} {2 n + 1} \cdot
  \binom{3n} n \cdot \frac{q_{2n + 1}(1)}{2n + 1} - q_1(1).$$
  This example is new.
\end{example}

\begin{example}
  Let $\gamma = \sqrt{\td}$ be the \emph{square root of the Todd
    class}, i.e.~$f = \sqrt{\frac x{1 - \exp(-x)}}$. Then
  $$g = 2 \Shi \left(\frac t 2\right) = \sum_{n = 0}^\infty \frac
  1{4^n \cdot (2n + 1) \cdot (2n + 1)!} \cdot t^{2n + 1},$$
  where
  $\Shi(y) := \int_0^y \frac{\sinh\eta}\eta d\eta$ is the hyperbolic
  sine integral.
  Therefore
  $$A(\log q(\sqrt{\td})) = \sum_{n \ge 0} \frac 1{4^n \cdot (2 n +
    1)!} \cdot \frac{q_{2n + 1}(1)}{2n + 1} - q_1(1).$$
  This example
  is new. (The calculation for $\td$ instead of $\sqrt{\td}$ leads to
  a much more complicated result.)
\end{example}

\begin{remark}
  As every product of Chern classes occurs in a suitable
  multiplicative class (e.g.~$\gamma = \prod_{i, j = 1}^\infty (1 +
  \lambda_i t_j) \in U_R$ with $R = \set Q[t_1, t_2, \dots]$), we can
  use the above corollary to calculate any characteristic class of the
  Hilbert schemes of points on the affine plane in terms of Nakajima's
  creation operators.
\end{remark}

We turn now from the tangent bundle to the multiplicative classes
built up from the $b_i$: Let $\gamma \in V_R$ be a multiplicative
characteristic class of the form $\gamma = \prod_{i = 1}^\infty
f(\lambda'_i)$ for a (uniquely determined) power series $f \in 1 + x
R[[x]]$.

We define a power series $g \in t R[[t]]$ by
$$
\frac {\partial g}{\partial t} \left(\frac x {f(- x)}\right) = f(- x),
$$
i.e.~$- t \frac {\partial g} {\partial t}$ is the compositional
inverse of the power series $- \frac x {f(x)}$.

\begin{theorem}
  Using the homomorphism $A$, the following partial result on
  $q(\gamma) \in Q_R$ holds, namely
  $$
  A(\log q(\gamma)) = \sum_{k \ge 1}^\infty \left(\frac 1 {k!}
    \left.\frac{\partial^k}{\partial t^k}\right|_{t = 0} g(t) \right)
  \cdot q_k(1) - q_1(1).$$
\end{theorem}

\begin{proof}
  We have
  $$A(q(\gamma))|_{(X, \alpha)} \one
  = q(\gamma)|_{(\set C^2, 1)} \one
  = \sum_{n \ge 0} \gamma|_{\Hilb n{(\set C^2, 1)}}.$$
  Thus again, everything boils down to calculating universal
  classes on the Hilbert schemes of points on the affine plane. As in
  the proof of the previous theorem, we have to show that
  $$\sum_{k \ge 1} g_k q_k(1)\vac = \sum_{n \ge 0} \gamma_{n -
    1}|_{\Hilb n {(\set C^2, 1)}}$$
  with $g = \sum_{k \ge 1} g_k t^k$.
  The right hand side of the equation is in fact a sum of
  characteristic classes of the tautological bundles $\Hilb n {\sheaf
    O}$ as $\ch(\Hilb n {\sheaf O}) = \Hilb n 1$ (note that
  $\ch(\sheaf O) = 1$). To be more precise, it is
  $$\gamma_{n - 1}|_{\Hilb n {(\set C^2, 1)}} = \left.\frac 1{(n -
      1)!}
    \frac{\partial^{n - 1}}{\partial x^{n - 1}}\right|_{x = 0}
  \prod_{i = 1}^n f(\Hilb n \epsilon_i \cdot x),$$
  where $\Hilb n \epsilon_1, \dots, \Hilb n \epsilon_n$ are the Chern
  roots of $\Hilb n {\sheaf O}$.

  From here, we continue as in the proof of the previous theorem.
  We just have to replace the
  tangent bundle by the tautological bundle. By~\cite{BOS1}, we have
  \begin{multline*}
    \gamma_{n - 1}|_{\Hilb n {(\set C^2, 1)}} = \sum_{\norm\lambda = n} \frac 1
    {h(\lambda)} \left(\left.\frac 1 {(n - 1)!} \frac{\partial^{n -
            1}}{\partial x^{n - 1}}\right|_{x = 0}
      \prod_{w \in
        D(\lambda)} f(h(w) x)\right)
    \\ \cdot \chi^\lambda_{(n)} \cdot z_{(n)}^{-1} \cdot
    q_n(1)\vac,
  \end{multline*}
  where the symbols have the same meaning as in the proof of the
  previous theorem. (Again, we have replaced the total Chern class by
  an arbitrary multiplicative class.) Doing the same calculations as
  above, one eventually arrives at
  $$\gamma_{n - 1}|_{\Hilb n {(\set C^2, 1)}} = \frac{q_n(1)} n \vac
  \cdot \left.\frac {(-1)^{n - 1}} {n!}
    \frac{\partial^{n - 1}}{\partial x^{n - 1}}\right|_{x = 0}
  (f(x))^n.$$
  Again by the Langrange Inversion Theorem this
  yields
  $$\sum_{n \ge 0} \gamma_{n - 1}|_{\Hilb n {(\set C^2, 1)}}
  = \sum_{k \ge 1} \left(\left.\frac{1}{k!}\frac{\partial^k}{\partial
        t^k}\right|_{t = 0} g\right) \cdot q_k(1) \vac,$$
  which again is all that remained to show.
\end{proof}

\begin{corollary}
  Consider the marked complex surface $(\set C^2, 1)$. We have
  $$
  \gamma|_{\Hilbtot (\set C^2, 1)} = \exp\left(
    \sum_{k \ge 1}^\infty \left(\frac 1 {k!}
      \left.\frac{\partial^k}{\partial t^k}\right|_{t = 0} g\right) \cdot
    q_k(1)\right)\vac.
  $$
\end{corollary}

\begin{remark}
  From this result, one can deduce a formula for
  $$\gamma|_{\Hilbtot (\mathbf C^2, r)} = \sum_{n \ge 0}
  \gamma|_{(\Hilb n X, \Hilb n r)}$$ for each $r \in \set Z$, as
  $\gamma(\Hilb n X, r^{[n]}) = \gamma(\Hilb n X, 1^{[n]})^r$. From
  this one can conclude that this formula holds for all $r \in R$.
\end{remark}

We do this in the following example:

\begin{example}
  Let $r \in \set Z$ and let $\gamma = c'^r$, i.e.~$f = (1 +
  x)^r$. Then
  $$\frac{\partial g}{\partial t}(x \cdot (1 - x)^{-r}) = (1 - x)^r.$$
  By the Lagrange Inversion Theorem,
  $$g = \sum_{n \ge 0} \frac{(-1)^{n - 1}}{n^2} \cdot \binom{r n}{n - 1}
  t^n.$$
  By the theorem, this gives
  $$A(\log q((c')^r)) = \sum_{n \ge 0} \frac{(-1)^{n - 1}}{n} \cdot
  \binom{r n}{n - 1} \cdot \frac{q_n(1)} n - q_1(1)$$ and by the
  corollary,
  $$c'|_{(\set C^2, r)} = \gamma|_{(\set C^2, 1)}
  = \exp\left(\sum_{n \ge 1} \frac{(-1)^{n - 1}}{n} \cdot
    \binom{r n}{n - 1} \cdot \frac{q_n(1)} n\right)\vac,$$
  which is an instance of Lehn's formula for $r = 1$,
  i.e.~corollary~\ref{cor:lehn_formula} applied to $\delta = 0$. .
  The case $r \ge 2$ appears to be new.
\end{example}

\begin{remark}
  One can use the previous corollary to obtain a complete description
  of the ring structure of the cohomology ring of the Hilbert schemes
  of points on the affine plane. (This has be studied successfully
  before, see~\cite{LS1}.)

  To understand this, consider the power series
  $$g = \sum_{k \ge 1} (\rho_k + \rho'_k) \cdot t^k \in t R[t]$$
  with
  $R := \set Q[\rho_1, \rho_2, \dots, \rho'_1, \rho'_2, \dots]$.
  There is a unique power series $f \in 1 + x R[[x]]$ such that
  $$\frac{\partial g}{\partial t}\left(\frac x {f(-x)}\right) = f(-x).$$
  Let $\gamma := \prod_{i = 1}^\infty f(\lambda'_i)$ be the
  multiplicative class defined by $f$. Then
  $$\gamma|_{\Hilbtot {(X, \alpha)}} = \exp\left(\sum_{k \ge 1}
    (\rho_k + \rho'_k) \cdot q_k(1)\right)\vac,$$
  i.e.~$\gamma|_{(X,
    \alpha)}$ can be seen as a \emph{universal cohomology class in
    $\set H_{\set C^2}$} as every class in $\set H_{\set C^2}$ can be
  written as a sum of polynomials in the $q_k(1)$ acting on the vacuum
  $\vac$.

  Now, $\gamma^2$ is another multiplicative class, which is defined by
  $f^2$ and which applied to $\Hilbtot{(X, \alpha)}$ gives
  $$\gamma^2|_{\Hilbtot{(X, \alpha)}} = \exp\left(\sum_{k \ge 1} (\rho_k +
    \rho'_k) \cdot q_k(1)\right)\vac \cdot \exp\left(\sum_{k \ge 1} (\rho_k +
    \rho'_k) \cdot q_k(1)\right)\vac,$$
  which includes all cup-products of all polynomials in the $q_k(1)$
  acting on the vacuum.

  Now let $h \in t R[[t]]$ be the unique power series with
  $$\frac{\partial h}{\partial t}\left(\frac x {f^2(-x)}\right) =
  f^2(-x).$$
  Then
  $$\left(\exp\left(\sum_{k \ge 1} (\rho_k + \rho'_k) \cdot
      q_k(1)\right)\vac\right)^2 = \exp\left(\sum_{k \ge 1}
    \left(\left.\frac 1{k!}  \frac{\partial^k}{\partial t^k}\right|_{t
        = 0} h(t)\right) \cdot q_k(1)\right)\vac.$$
  Thus, we can
  effectively calculate all powers of the universal cohomology class,
  which encodes completely the ring structure.
\end{remark}

\appendix

\section{Lagrange Inversion Theorem}

For the convenience of the reader, we state the Langrange Inversion
Theorem (also called Lagrange expansion) in the form in which it is
needed in the main text. For a proof of the general statement see
e.g.~\cite{WW}.

Let $f \in 1 + x R[[x]]$ be a power series and define $g \in t R[[t]]$
such that
$$\frac{\partial g}{\partial t}\left(\frac x f\right) = f,$$
i.e.~such that $t \frac{\partial g}{\partial t}$ is the compositional
inverse of $\frac x f$.
\begin{proposition}
  \label{prop:lagrange_expansion}
  The power series $g$ is given by
  $$g = \sum_{n = 1}^\infty \frac{t^n}{n \cdot n!}
  \left(\left.\frac{\partial^{n -
          1}}{\partial x^{n - 1}}\right|_{x = 0} f^n\right).$$
\end{proposition}
\qed

\section{Two algebraic identities}

The following identity holds true in the domain of the rationals for
all $m \in \set N_0$ and $0 \leq p \leq m$:
\begin{lemma}
  It is $$\sum_{s = 0}^m \frac{(-1)^s s^p}{s! (m-s)!} = \begin{cases}
    0 & \text{for $p < m$}
    \\
    (-1)^m & \text{for $p = m$}.
  \end{cases}
  $$
\end{lemma}

\begin{proof}
  Observe that
  \begin{align*}
    \sum_{s = 0}^m \frac{(-1)^s s^p}{s! (m - s)!}
    & = \left.\sum_{s = 0}^m \frac{(-x)^s s^p}{s! (m - s)!}\right|_{x
      = 1}
    \\
    & = \left.\left(x \frac d{dx}\right)^p \sum_{s = 0}^m
      \frac{(-x)^s}{s!(m - s)!}\right|_{x = 1}
    \\
    & = \frac 1 {m!} \left.\left(x \frac d{dx}\right)^p
      (1-x)^m\right|_{x = 1},
  \end{align*}
  where $x$ is a formal variable. As
  $$\left(x \frac d {dx}\right)^p (1-x)^m$$
  is divisable as a
  polynomial in $x$ by $(1 - x)$ for $p < m$ and $$\left(x \frac d {d
      x}\right)^m (1 - x)^m = m! \cdot x^n + f$$
  where $f$ is a polynomial
  divisable by $1 - x$, the lemma follows.
\end{proof}

Let $f \in 1 + x\set R[[x]]$ any formal power series with constant
coefficient one. For each $n \in \set N_0$ define the formal power
series 
$$P_n := \sum_{s = 0}^n \frac{(-1)^s}{s!(n - s)!} \cdot \prod_{k = -(n -
  s)}^s f(k \cdot x)$$
in $1 + x R[[x]]$.
\begin{lemma}
  \label{lem:identity2}
  The power series $P_n$ has leading term $(-1)^n \left(\left.\frac 1 {n!}
      \frac {d}{d x^n}\right|_{x = 0} f^{n + 1}\right) \cdot
  x^n$, i.e.~
  $$P_n = (-1)^n \cdot \left(\left.\frac {d}{d x^n}\right|_{x
      = 0} f^{n + 1}\right) \cdot x^n + O(x^{n + 1}).$$
\end{lemma}

\begin{proof}
  First of all, we consider the logarithm
  $$\log f = \sum_{i = 1}^\infty F_i x^i$$
  with certain coefficients
  $F_i \in R$ as a formal power series in $xR[[x]]$.
  This enables us to write
  \begin{align*}
    \prod_{k = - (n - s)}^s f(k \cdot x)
    & = \exp\left(\sum_{k = -(n - s)}^s (\log f)(k \cdot x)\right)
    \\
    & = \exp\left(\sum_{i = 1}^\infty \sum_{k = -(n - s)}^s F_i k^i
      x^i\right)
    \\
    & = \exp\left(\sum_{i = 1}^\infty F_i x^i \left(\sum_{k = 1}^s k^i
        + \sum_{k = 1}^{n - s} (-k)^i\right)\right).
  \end{align*}
  Recall the classical summation formula involving the Bernoulli
  numbers. These show, that terms like $\sum_{k = 1}^s k^i$ can be
  expressed as polynomials in $s$ of degree $i + 1$, in fact
  $$\sum_{k = 1}^s k^i = \frac 1 {i + 1} \left(s^{i + 1} + \frac 1 2
    (i + 1) s^i + \dots\right),$$
  where ``$\dots$'' denotes term of lower degree in $s$. Similarly,
  one has
  \begin{multline*}
    \sum_{k = 1}^{n - s} (-k)^i
    = \frac{(-1)^i} {i + 1}\left((n -
      s)^{i + 1}
      + \frac 1 2 (i + 1) (n - s)^i + \dots\right)
    \\
    = \frac 1{i + 1}\left(-s^{i + 1} + \frac 1 2 (i + 1)(2n + 1) s^i +
      \dots\right).
  \end{multline*}
  Summing up yields
  $$\sum_{k = 1}^s k^i + \sum_{k = 1}^{n - s} (-k)^i 
  = (n + 1) s^i + \dots.$$
  Thus, we have
  $$\prod_{k = - (n - s)}^s f(k \cdot x) = \exp\left(\sum_{i =
        1}^\infty ((n + 1) F_i x^i (s^i + \dots))\right)
    \\
    = f^{n + 1}(s \cdot x) + g(s, x)
  $$
  for a power series $g \in R[s][[x]]$, in which we collect all the
  ``$\dots$''-terms. In each homogeneous component of $g$ the exponent
  of $s$ is strictly less than the exponent of $x$.
  Plugging this into our definition of $P_n$, we arrive at
  $$P_n = \sum_{s = 0}^n \frac{(-1)^s} {s!(n - s)!} f^{n + 1}(s \cdot
  x) + x g(s \cdot x).$$
  In order to complete the proof, we have to
  differentiate this expression up to $n$-times with respect to
  $x$. So, let $0 \leq p \leq n$. Then
  $$
  \left.\frac 1 {p!} \frac d{dx^p}\right|_{x = 0} P_n = \sum_{s =
    0}^n \frac{(-1)^s}{s!(n - s)!} \left. \frac 1 {p!} \frac d
    {dx^p}\right|_{x = 0}(f^{n + 1}(s \cdot x) + x g(s \cdot x)).$$
  The term $\left.\frac 1 {p!} \frac{d}{dx^p}\right|_{x = 0} f^{n +
    1}(s \cdot x)$ is of degree $p$ (or less) in $s$, the term
  $\left.\frac 1 {p!} \frac{d}{dx^p}\right|_{x = 0} x g(s \cdot x)$ is
  of degree strictly less than $p$ in $s$. Therefore, in view of the
  previous lemma, we have
  $$
  \left.\frac 1 {p!} \frac d{dx^p}\right|_{x = 0} P_n = 
  \begin{cases}
    (-1)^n \left.\frac{1}{n!}\frac{d}{d x^n}\right|_{x = 0} f^{n + 1}
    & \text{for $p = n$}
    \\
    0 & \text{for $p < n$}.
  \end{cases}$$
\end{proof}

\bibliographystyle{amsplain}
\bibliography{universal_classes}

\end{document}